\documentclass[11pt]{article}

\newskip\stdskip                      
\usepackage{amscd}     
\usepackage{amsfonts}
\usepackage{amsmath}
\usepackage{amssymb}
\usepackage{amstext}
\usepackage{amsthm}          
\usepackage{graphicx}
\usepackage{psfrag}
\usepackage[all]{xy}
\usepackage{url}

\newcommand{\R}{\mathbb R}

\newcommand{\Q}{\mathbb Q}
\newcommand{\Z}{\mathbb Z} 
\newcommand{\F}{\mathcal F}
\newcommand{\rk}{{\rm rk} \,}
\newcommand{\spc}{{\rm Spin}^c}
\newcommand{\gt}{\mathfrak{t}}

\newtheorem{thm}{Theorem}[section]
\newtheorem{lemma}[thm]{Lemma}
\newtheorem{cor}[thm]{Corollary}
\newtheorem{prop}[thm]{Proposition}
\newtheorem{conj}[thm]{Conjecture}

\theoremstyle{definition}
\newtheorem{dfn}[thm]{Definition}
\newtheorem{rem}[thm]{Remark}

\begin{document}
\author{Paolo Ghiggini \\ 
CIRGET, Universit\'e du Qu\'ebec \`a Montr\'eal \\
Case Postale 8888, succursale Centre-Ville \\
Montr\'eal (Qu\'ebec) H3C 3P8, Canada \\
\url{ghiggini@math.uqam.ca}}
\title{Knot Floer homology detects genus-one fibred knots}
\date{\today}
\maketitle

\begin{abstract} \noindent
Ozsv\'ath and Szab\'o conjectured that knot Floer homology 
detects fibred knots. We propose a strategy to approach this 
conjecture based on Gabai's theory of sutured manifold 
decomposition and contact topology. We implement this 
strategy for genus-one knots, obtaining as a corollary 
that if rational surgery on a knot $K$ gives the Poincar\'e homology 
sphere $\Sigma(2,3,5)$, then $K$ is the left-handed trefoil knot.
\end{abstract}

\section{Introduction}
Knot Floer homology was introduced independently by Ozsv\'ath and Szab\'o
in \cite{O-Sz:knots} and by Rasmussen in \cite{rasmussen:tesi}.
For any knot $K$ in $S^3$ and any integer $d$ the knot Floer homology 
group $\widehat{HFK}(K, d)$ is a finitely generated graded 
Abelian group.  

Knot Floer homology can be seen as a categorification of the Alexander 
polynomial in the sense that 
\[ \sum \chi(\widehat{HFK}(K,d))T^i = \Delta_K(T) \]
where $\Delta_K(T)$ denotes the symmetrised Alexander polynomial. However
the groups $\widehat{HFK}(K,d)$, and in particular the top non trivial
group, contain more information than just
the Alexander polynomial, as the following results show. 
\begin{thm} (\cite[Theorem 1.2]{O-Sz:genus})\label{genere}
Let $g(K)$ denote the genus of $K$. Then
\[ g(K)= \max \big \{ d \in \Z : \widehat{HFK}(K, d) \neq 0 \big \}. \]  
\end{thm}
\begin{thm} (\cite[Theorem 1.1]{O-Sz:cont}) \label{fibrato}
Let $g$ be the genus of $K$. If $K$ is a fibred knot, then 
$\widehat{HFK}(K, g)= \Z$.
\end{thm}
 
Ozsv\'ath and Szab\'o formulated the following conjecture, whose evidence
is supported by the computation of knot Floer homology for a large 
number of knots. 

\begin{conj}\label{congettura} (Ozsv\'ath--Szab\'o)
If $K$ is a knot in $S^3$ with genus $g$ and $\widehat{HFK}(K, g)= \Z$
then $K$ is a fibred knot.
\end{conj}

In this article we propose a strategy to attack Conjecture 
\ref{congettura}, and we implement it in the case of genus-one knots.  
More precisely, we will prove the 
following result.

\begin{thm}\label{principale}
Let $K$ be an oriented genus-one knot in $S^3$. 
Then $K$ is fibred if and only if $\widehat{HFK}(K, 1)= \Z$
\end{thm}

Our strategy to prove Theorem \ref{principale} is to deduce 
information about the top knot Floer homology group
of non-fibred knots from topological properties of their 
complement via sutured manifold decomposition and
contact structures in a way that is reminiscent of the proof of 
Theorem \ref{genere}. 

It is well known that the only fibred knots of genus one are the 
trefoil knots and the figure-eight knot, therefore Theorem 
\ref{principale}, together with the computation of knot Floer homology
for such knots, implies the following
\begin{cor}\label{nicecorollary1}
Knot Floer homology detects the trefoil knots and the 
figure-eight knot.
\end{cor}
The following conjecture was formulated by Kirby in a remark after 
Problem 3.6(D) of his problem list, and by Zhang in \cite{zhang:I}.
\begin{conj}(Conjecture $\hat{I}$)
If $K$ is a knot in $S^3$ such that there exists a rational 
number $r$ for which the $3$--manifold obtained by 
$r$--surgery on $K$ is homeomorphic to the Poincar\'e homology sphere
$\Sigma(2,3,5)$, then $K$ is the left-handed trefoil knot.
\end{conj}
Conjecture $\hat{I}$ was proved for some knots by Zhang in 
\cite{zhang:I}, and a major step toward its complete proof was made 
by Ozsv\'ath and Szab\'o in 
\cite{O-Sz:seifert}, where they proved that a counterexample to 
Conjecture $\hat{I}$ must have the same knot Floer Homology groups as
the left-handed trefoil knot.
Corollary \ref{nicecorollary1} provides the missing step to prove it 
in full generality.
\begin{cor}\label{nicecorollary2}
Conjecture $\hat{I}$ holds.
\end{cor}

Corollary \ref{nicecorollary1} has also been used by Ozsv\'ath and Szab\'o
to prove that the trefoil knot and the figure-eight knot are determined
by their Dehn surgeries \cite{O-Sz:trefoil}.

\subsection*{Acknowledgements}
We warmly thank Steve Boyer, Ko Honda, Joseph Maher, and Stefan 
Tillmann for many inspiring conversations. 

\section{Overview of Heegaard Floer theory}
Heegaard Floer theory is a family of invariants introduced by Ozsv\'ath
and Szab\'o in the last few years for the most common
objects in low-dimensional topology. In this section we will give a 
brief overview of the results in Heegaard Floer theory we will need 
in the following, with no pretension of completeness. The details can 
be found in Ozsv\'ath and Szab\'o's papers 
\cite{O-Sz:1,O-Sz:2,O-Sz:3,O-Sz:knots,O-Sz:cont,O-Sz:genus}.

\subsection{Heegaard Floer homology}
Let $Y$ be a closed, connected, oriented $3$--manifold. For any
$\spc$--structure $\gt$ on $Y$ Ozsv\'ath and Szab\'o 
\cite{O-Sz:1} defined an Abelian group $HF^+(Y, \gt)$ which is
an isomorphism invariant of the pair $(Y, \gt)$. When $c_1(\gt)$ is 
not a torsion element in $H^2(Y)$ the group $HF^+(Y, \gt)$ is finitely
generated; see \cite[Theorem 5.2 and Theorem 5.11]{O-Sz:2}. 
If there is a distinguished 
surface $\Sigma$ in $Y$ which is clear from the context, we use the 
shortened notation 
\[ HF^+(Y,d)= \bigoplus_{\begin{array}{c} \gt \in \spc(Y) \\ \langle c_1(\gt), [\Sigma] \rangle =2d
\end{array}} HF^+(Y, \gt). \]
This notation makes sense because $HF^+(Y, \gt) \neq 0$ only for finitely
many $\spc$--structures. Heegaard Floer homology is symmetric; in fact
there is a natural isomorphism $HF^+(Y, d) \cong HF^+(Y, -d)$ for any 
$3$--manifolds $Y$ and any integer $d$: see \cite[Theorem 2.4]{O-Sz:2}. 
There is an adjunction inequality relating $HF^+$ to the minimal genus 
of embedded surfaces which can be stated as follows.

\begin{thm} \cite[Theorem 1.6]{O-Sz:2} \label{aggiunzione}
If $\Sigma$ has genus $g$, then $HF^+(Y,d)=0$ for all $d > g$.
\end{thm}

Any connected oriented cobordism $X$ from $Y_1$ to 
$Y_2$ induces a homomorphism 
\[ F_X \colon HF^+(Y_1, \gt_1) \to HF^+(Y_2, \gt_2) \]
which splits as a sum of homomorphisms indexed by the 
$\spc$--structures on $X$ extending $\gt_1$ and $\gt_2$. When
$X$ is obtained by a single $2$--handle addition $F_X$ fit into
a surgery exact triangle as follows.
 
\begin{thm} (\cite[Theorem 9.12]{O-Sz:2}) \label{triangolo}
Let $(K, \lambda)$ be an oriented framed knot in $Y$, and let $\mu$ a 
meridian of $K$. Denote by $Y_{\lambda}(K)$ the manifold obtained from 
$Y$ by surgery on $K$ with framing $\lambda$, and by $X$ the cobordism
induced by the surgery. If we choose a surface $\Sigma \subset Y \setminus K$ to
partition the $\spc$-structures, then the following triangle \\
\begin{center}
\parbox{10cm}{\xymatrix{ 
HF^+(Y,d) \ar^{F_X}[rr] & & HF^+(Y_{\lambda}(K), d) \ar[dl] \\
  & HF^+(Y_{\lambda + \mu }(K), d) \ar[ul]
}} \end{center} 
is exact for any $d \in \Z$.
\end{thm}

\subsection{The Ozsv\'ath--Szab\'o contact invariant}
A contact structure $\xi$ on a $3$--manifold $Y$ determines a 
$Spin^c$--structure $\mathfrak{t}_{\xi}$ on $Y$ such that 
$c_1(\mathfrak{t}_{\xi})=c_1(\xi)$. To any contact manifold $(Y, \xi)$ we can 
associate an element $c^+(\xi) \in HF^+(-Y, \mathfrak{t}_{\xi})/ \pm 1$ 
which is an isotopy invariant of $\xi$, see \cite{O-Sz:cont}.
In the following we will always abuse the notation and consider 
$c^+(\xi)$ as an element of $HF^+(-Y, \mathfrak{t}_{\xi})$, 
although it is, strictly speaking, defined only up to sign. This
abuse does not lead to mistakes as long as we do not use the 
additive structure on $HF^+(-Y, \mathfrak{t}_{\xi})$.

The proof of the following lemma is contained in the proof of
\cite[Corollary~1.2]{O-Sz:genus}; see also 
\cite[Theorem~2.1]{kmosz} for a similar result in the setting
 of monopole Floer homology.

\begin{lemma} \label{nonzero}
Let $Y$ be a closed, connected oriented $3$--manifold with
$b_1(Y)=1$, and let $\xi$ be a weakly symplectically fillable
contact structure on $Y$ such that $c_1(\xi)$ is non trivial in 
$H^2(Y; \R)$. Then 
$c^+(\xi)$ is a primitive element of $HF^+(-Y, \gt_{\xi})/ \pm 1$.
\end{lemma}

Given a contact manifold $(Y, \xi)$ and a Legendrian knot $K \subset Y$
there is an operation called {\em contact $(+1)$--surgery} which produces
a new contact manifold $(Y', \xi')$; see \cite{ding-geiges:1} and 
\cite{ding-geiges:2}. The Ozsv\'ath--Szab\'o contact invariant behaves 
well with respect to contact $(+1)$--surgeries.

\begin{lemma} (\cite[Theorem 2.3]{lisca-stipsicz:3}; see also
\cite{O-Sz:cont}) \label{surgery}
Suppose $(Y', \xi')$ is obtained from $(Y, \xi)$ by a contact 
$(+1)$--surgery. Let $-X$ be the cobordism induced by the surgery with
opposite orientation. Then
\[ F_{-X}^+ (c^+(\xi)) = c^+(\xi'). \]
\end{lemma}

\subsection{Knot Floer homology}
Knot Floer homology is a family of finitely dimensional graded Abelian 
groups $\widehat{HFK}(K, d)$ indexed by $d \in \Z$ attached to any 
oriented knot $K$ in $S^3$; see Ozsv\'ath and Szab\'o \cite{O-Sz:knots}. 
Denote the $0$--surgery on $K$ by $Y_K$. The knot Floer homology of $K$
is related to the Heegaard Floer homology of $Y_K$ by the following
 
\begin{prop} \label{relazione} (\cite[Corollary 4.5]{O-Sz:knots} and 
 \cite[Corollary 1.2]{O-Sz:genus}).
Let $K$ be a knot of genus $g>1$. Then
\[ \widehat{HFK}(K, g)= HF^+(Y_K, g-1). \]
\end{prop}

Another property of knot Floer homology we will need later is a 
K\"unneth-like formula for connected sums.
\begin{prop} \label{kunneth} (\cite[Corollary~7.2]{O-Sz:knots})
Let $K_1$ and $K_2$ be knots in $S^3$, and denote by $K_1 \# K_2$ their 
connected sum. If $\widehat{HFK}(K_1, d)$ is a free Abelian group for
every $d$ (or if we work with coefficients in a field), then 
\[  \widehat{HFK}(K_1 \# K_2, d)= \bigoplus_{d_1+d_2=d} \widehat{HFK}(K_1, d_1) \otimes
\widehat{HFK}(K_1, d_1). \]
\end{prop} 

The definition of Knot Floer homology can be extended to links.
To any link $L$ in $S^3$ with $|L|$ components, Ozsv\'ath and Szab\'o 
associate a null-homologous knot $\kappa(L)$ in $\#^{|L|-1} (S^2 \times S^1)$. For a link
with two components the construction is the following: choose points
$p$ and $q$ on different components of $L$, then replace two balls 
centred at $p$ and $q$ with a $1$--handle $S^2 \times [0,1]$ and define 
$\kappa(L)$ as the banded connected sum of the two components of $L$ 
performed with a standardly embedded band in the $1$--handle. If $L$ 
has more connected components this operation produces a link with one 
component less, so we repeat it until we obtain a knot. For the 
details see \cite[Section 2.1]{O-Sz:knots}.

\begin{dfn}(\cite[Definition 3.3]{O-Sz:knots})
If $L$ is a link in $S^3$ we define 
\[ \widehat{HFK}(L,d)= \widehat{HFK}(\kappa(L), d). \]
\end{dfn}

\section{Taut foliations and Heegaard Floer 
homology}
\subsection{Controlled perturbation of taut foliations}
Eliashberg and Thurston in \cite{eliashberg-thurston} introduced a new
technique to construct symplectically fillable contact structures by 
perturbing taut foliations. In this section we show how to control the
perturbation in the neighbourhood of some closed curves. 
We need to introduce some terminology about confoliations, 
following Eliashberg and Thurston \cite{eliashberg-thurston}.

\begin{dfn} 
A {\em confoliation} on an oriented $3$--manifold
is a tangent plane field $\eta$ defined by a $1$--form $\alpha$ such that
$\alpha \land d \alpha \geq 0$.
\end{dfn}
Given a confoliation $\eta$ on $M$ we define its 
{\em contact part} $H(\eta)$ as 
\[ H(\eta)= \{ x \in M : \alpha \land d \alpha (x) > 0 \}. \]

\begin{lemma}\label{olonomia}
Let $\Sigma$ be a compact leaf with trivial germinal holonomy in a taut 
smooth foliation $\F$ on a $3$--manifold $M$, and let $\gamma$ be a 
non-separating closed curve in $\Sigma$. Then we can modify $\F$ in a 
neighbourhood of $\gamma$ so that we obtain a new taut smooth foliation 
with non trivial linear holonomy along $\gamma$.
\end{lemma}
\begin{proof}
The holonomy of $\Sigma$ determines the germ of $\F$ along $\Sigma$ (see
\cite[Theorem 3.1.6]{candel-conlon:1}), therefore $\Sigma$ has a 
neighbourhood $N = \Sigma \times [-1,1]$ such that $\F|_{N}$ is the product 
foliation. Pick $\gamma' \subset \Sigma = \Sigma \times \{ 0 \}$ such that it intersects $\gamma$
in a unique point, and call $\Sigma' = \Sigma \setminus \gamma'$. The boundary of $\Sigma'$ has
two components $\gamma'_+$ and $\gamma_-'$. For every point $x \in \gamma'$ denote by 
$x_{\pm}$ the points in $\gamma_{\pm}'$ which correspond to $x$. 
Choose a diffeomorphism
$f \colon [-1,1] \to [-1,1]$ such that
\begin{enumerate}
\item $f(x)=x$ if $x \in [-1, -1+ \epsilon] \cup [1- \epsilon, 1]$ for some small $\epsilon$;
\item $f(0)=0$;
\item $f'(0) \neq 1$,
\end{enumerate} 
then re-glue $\gamma_-' \times [-1,1]$ to $\gamma_+' \times [-1,1]$ identifying $(x_-, t)$
to $(x_+, f(t))$.
\end{proof}

\begin{lemma}\label{pippo}
Let $(M, \F)$ be a foliated manifold, and let $\gamma$ be a curve
with non-trivial linear holonomy contained in a leaf $\Sigma$.
Then $\F$ can be approximated in the $C^0$--topology by  
 confoliations such that $\gamma$ is contained in their contact 
parts and is a Legendrian curve with twisted number zero 
with respect to the framing induced by $\Sigma$. 
\end{lemma}

\begin{proof}
We apply \cite[Proposition~2.6.1]{eliashberg-thurston} to
make $\F$ a contact structure in a neighbourhood of $\gamma$, then 
the proof of 
the lemma is a check on the explicitly given contact form.
\end{proof}

For any subset $A \subset M$ we define its {\em saturation} 
$\widehat{A}$ as the set of all points in $M$ which can
be connected to $A$ by a curve tangent to $\eta$. 

\begin{dfn}
A confoliation $\eta$ on $M$ is called {\em transitive} if 
$\widehat{H(\eta)}=M$.
\end{dfn}

\begin{lemma} \label{controverso}
Let $(M, \F)$ be a smooth taut foliated manifold, let $\Sigma$ be a 
compact leaf of $\F$ with trivial germinal holonomy, and let $\gamma \subset \Sigma$ 
be a closed non-separating curve. Then $\F$ can be approximated in the 
$C^0$--topology by contact structures such that $\gamma$ is a Legendrian 
curve with twisting number zero with respect to $\Sigma$.
\end{lemma}
\begin{proof}
First we apply Lemma \ref{olonomia} to create non-trivial 
linear holonomy along $\gamma$, so that we can apply Lemma \ref{pippo}
 to make $\F$ a contact structure in a neighbourhood of $\gamma$, and $\gamma$ 
becomes a Legendrian curve with twisting number zero.

The approximation of a confoliation $\F$ by contact structures is 
done in two steps. First $\F$ is $C^0$--approximated by a 
transitive confoliation $\widetilde{\F}$, then $\widetilde{\F}$ 
is $C^1$--approximated by a contact structure $\xi$.
The first step is done by perturbing the confoliation in 
arbitrarily small neighbourhoods of curves contained in 
$M \setminus \widehat{H(\F)}$, then we can assume that a neighbourhood $V$
of $\gamma$ is not touched in the first step. 

Since $\widetilde{\F}$ and 
$\xi$ are $C^1$--close, they are defined by $C^1$--close $1$--forms 
$\alpha$ and $\beta$. Let $h \colon M \to [0,1]$ be a smooth function supported
in $V$ such that $h \equiv 1$ in a smaller neighbourhood of $\gamma$, then 
the $1$--form $\beta + h(\alpha - \beta)$ coincides with $\alpha$ near $\gamma$ and with 
$\beta$ outside $V$. It defines a contact structure which is 
$C^0$--close to $\F$ because $\alpha$ and $\beta$ are 
$C^1$--close and the contact condition open in the $C^1$--topology. 
\end{proof}

\begin{prop}\label{fill} (\cite[Corollary~3.2.5]{eliashberg-thurston} 
and \cite[Corollary~1.4]{eliashberg-fill})
Let $\F$ be a taut foliation on $M$. Then any contact structure $\xi$
which is sufficiently close to $\F$ as a plane field in the $C^0$
topology is weakly symplectically fillable.
\end{prop} 

\subsection{An estimate on the rank of $HF^+(Y)$ coming from taut 
foliations} 
Let $\eta$ be a field of tangent planes in the $3$--manifold $Y$, 
and let $S$ be an embedded compact surface with non empty 
boundary such that $\partial S$ is tangent to $\eta$. 
 
\begin{dfn}
Let $v$ be the positively oriented unit tangent vector field of $\partial S$.
 We define the {\em relative Euler class} $e(\eta, S)$ of $\eta$ on $S$ as the
obstruction to extending $v$ to a nowhere vanishing section of $\eta|_S$.
\end{dfn}
Properly speaking $e(\eta, S)$ is an element of $H^2(S, \partial S)$, but we can 
identify it with an integer number via the isomorphism 
$H^2(S, \partial S) = \Z$.
If $\eta$ is the field of the tangent planes of the leaves of a 
foliation $\F$ we write $e(\F, S)$ for $e(\eta, S)$. If $\eta$ is a 
contact structure $e(\eta, S)$ is the sum of the rotation numbers of 
the components of $\partial S$ computed with respect to $S$.

\begin{figure}\centering
\psfrag{+}{\footnotesize $\alpha_+$}
\psfrag{-}{\footnotesize $\alpha_-$}
\psfrag{S}{\footnotesize $S$}
\includegraphics[width=5cm]{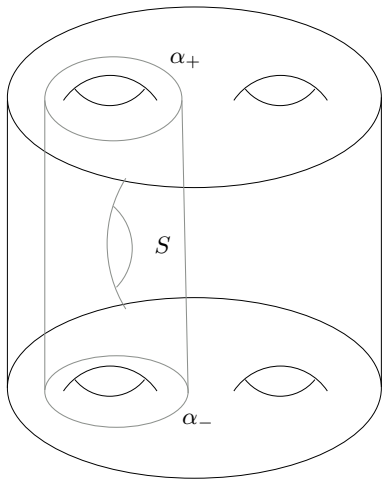}
\caption{A somehow misleading picture of $S \subset Y \setminus \Sigma$}
\label{S.fig}
\end{figure}

\begin{thm}\label{nonnapapera}
Let $Y$ be a $3$--manifold with $H_2(Y) \cong \Z$, and let $\Sigma$ be a genus 
minimising closed surface representing a generator of $H_2(Y)$. 
Call $\Sigma_+$ and $\Sigma_-$ the two components of $\partial (Y \setminus \Sigma)$. Suppose that
$\Sigma$ has genus $g(\Sigma)>1$ and that $Y$ admits two 
smooth taut foliations $\F_1$ and $\F_2$ such that $\Sigma$ is a compact 
leaf for both, and the holonomy of $\Sigma$ has the same Taylor series as 
the identity. If there exists a properly embedded surface $S \subset Y \setminus \Sigma$ 
with boundary $\partial S= \alpha_+ \cup \alpha_-$ such that 
\begin{enumerate}
\item $\alpha_+ \subset \Sigma_+$ and $\alpha_- \subset \Sigma_-$ are non separating curves, and
\item $e(\F_1, S) \neq e(\F_1, S)$
\end{enumerate}
(see Figure \ref{S.fig}), then $\rk HF^+(Y, g-1) > 1$.
\end{thm} 

Requiring that the holonomy of $\Sigma$ has the same Taylor series as 
the identity is not as strong a restriction as it seems, because 
foliations constructed by sutured manifold theory have this property;
see the induction hypothesis in the proof of 
\cite[Theorem~5.1]{gabai:1}.

The strategy of the proof is to view $\F_1$ and $\F_2$ as taut foliations
on $-Y$ and to approximate them by contact structures $\xi_1$ and $\xi_2$ on 
$-Y$ so that $c(\xi_i) \in HF^+(Y, g-1)$, then to construct a new 
$3$--manifold $Y_{\phi}$ together with a $4$--dimensional cobordism $W$ 
from $-Y$ to $-Y_{\phi}$ so that $F_{-W}^+(c^+(\xi_1))$ and $F_{-W}^+(c^+(\xi_2))$ are 
linearly independent in $HF^+(Y_{\phi})$. This implies that $c^+(\xi_1)$ and 
$c^+(\xi_2)$ are linearly independent in $HF^+(Y, g-1)$.
The entire subsection is devoted to the proof of Theorem 
\ref{nonnapapera}. 

We choose a diffeomorphism $\phi \colon \Sigma_+ \to \Sigma_-$ such that $\phi(\alpha_+)= \alpha_-$, and 
we form the new $3$--manifold $Y_{\phi}$ by cutting $Y$ along $\Sigma$ and 
re-gluing $\Sigma_+$ to $\Sigma_-$ after acting by $\phi$. The diffeomorphism $\phi$
exists because $\alpha_+$ and $\alpha_-$ are non separating. 

It is well known that the mapping class group of a closed surface
of genus $g>1$ is generated, as a monoid, by positive Dehn twists;
see \cite[Footnote 3]{eliashberg-fill}. In order to 
construct the $4$--dimensional cobordism $W$ from $Y$ to $Y_{\phi}$ we 
decompose $\phi$ as a product $\phi= \prod \tau_{c_1} \ldots \tau_{c_k}$ where $\tau_{c_i}$ is 
a positive Dehn twist around a curve $c_i \subset \Sigma$. Then we identify
a tubular neighbourhood $N$ of $\Sigma$ with $\Sigma \times [-1,1]$, 
we choose distinct points $t_1, \ldots ,t_k$ in $(-1,1)$ and 
see $c_i$ as a curve in $\Sigma \times \{ t_i \}$. 
The surface $\Sigma \times \{ t_i \}$ induces  a framing on $c_i$, and 
$Y_{\phi}$ is obtained by $(-1)$--surgery on the link 
$C= c_1 \cup \ldots \cup c_k \subset Y$, where the surgery coefficient of $c_i$ is 
computed with respect to the framing induced by $\Sigma \times \{ t_i \}$.
 Equivalently, $- Y_{\phi}$ is obtained by 
$(+1)$--surgery on the same link $C$ seen as a link in $-Y$. 
We denote by $W$ the smooth $4$--dimensional cobordism obtained 
by adding $2$--handles to $-Y$ along the curves $c_i$ with 
framing $+1$, and by $-W$ the same cobordism with opposite 
orientation, so that $-W$ is obtained by adding $2$--handles 
to $Y$ along the curves $c_i$ with framing $-1$.
  
\begin{lemma} \label{isom}
The map 
\begin{equation*} 
 F^+_{-W} \colon HF^+(Y, g-1) \to HF^+(Y_{\phi}, g-1) 
\end{equation*}
induced by the cobordism $-W$ is an isomorphism.
\end{lemma}
\begin{proof} The map $F_{-W}^+$ is a composition of maps induced by
elementary cobordisms obtained from a single $2$--handle addition. 
We apply the surgery exact triangle (Theorem \ref{triangolo}) 
and the adjunction inequality (Theorem \ref{aggiunzione}) as in 
\cite[Lemma 5.4]{O-Sz:symp} to prove that each elementary cobordism 
induces an isomorphism in Heegaard Floer homology.
\end{proof}

We can assume that the tubular neighbourhood $N$ of $\Sigma$ is foliated 
as a product in both foliations. In fact there is a diffeomorphism 
$Y \cong (Y \setminus \Sigma) \cup (\Sigma \times [-1,1])$ such that $\Sigma_+$ is identified with
$\Sigma \times \{ -1 \}$, and $\Sigma_-$ is identified with $\Sigma \times \{ 1 \}$,  then we can 
extend the foliations $\F_i|_{Y \setminus \Sigma}$ to foliations on $(Y \setminus {\Sigma}) \cup
(\Sigma \times [-1,1])$ which are product foliations on $\Sigma \times [-1,1]$.
 We call the resulting  foliated manifolds
$(Y, \F_1)$ and $(Y, \F_2)$ again. 
This operation does not destroy the smoothness of $\F_1$ and $\F_2$ 
because their holonomies along $\Sigma$  have the same Taylor 
series as the identity.

We see $S$ as a surface in $Y \setminus N$, so that $\alpha_+$ is identified to a 
curve in $\Sigma \times \{ -1 \}$, and $\alpha_-$ is identified to a curve in 
$\Sigma \times \{ 1 \}$. By Lemma \ref{controverso} we can control the 
perturbations of $\F_1$ and $\F_2$ so that $\alpha_+$, $\alpha_-$, and $c_i$ for all 
$i$ become Legendrian curves with twisting number zero for both 
$\xi_1$ and $\xi_2$, where the twisting number is computed with respect 
to the framing induced by $\Sigma$. This implies that we can construct 
contact structures $\xi_1'$ and $\xi_2'$ on $-Y_{\phi}$ by $(+1)$--contact
 surgery on $\xi_1$ and $\xi_2$.

By hypothesis $b_1(Y)=1$, and $c_1(\xi_1)$ and $c_1(\xi_2)$ are non torsion 
because 
\[ \langle c_1(\xi_i), [\Sigma] \rangle = \langle c_1(\F_i), [\Sigma] \rangle = \chi(\Sigma) <0, \]
therefore Lemma \ref{nonzero} applies and gives $c^+(\xi_i) \neq 0$ 
for $i=1,2$. Moreover $c^+(\xi_1')= F^+_{-W}(c^+(\xi_1))$ and
$c^+(\xi_1')= F^+_{-W}(c^+(\xi_1))$ by Lemma \ref{surgery} because 
$\xi_i'$ is obtained form $\xi_i$ by a sequence of contact $(+1)$--surgeries,
therefore from Lemma \ref{isom} it follows $c^+(\xi_i') \neq 0$ for
$i=1,2$.

\begin{lemma}\label{diverse}
Let $N_{\phi}$ be the subset of $Y_{\phi}$ diffeomorphic to 
$\Sigma \times [-1,1]$ obtained from $N \subset Y$ by performing 
$(-1)$--surgery on $C=c_1 \cup \ldots \cup c_n$, and call $S'$ the annulus
in $N_{\phi}$ bounded by the curves $\alpha_+$ and $\alpha_-$. If we define 
$\overline{S}= S \cup S'$, then 
\[ \langle c_1(\xi_1'), [\overline{S}] \rangle \neq \langle c_1(\xi_2'), [\overline{S}] \rangle. \]
\end{lemma}
 
\begin{proof}
Because $\alpha_+$ and $\alpha_-$ are Legendrian curves with twisting number $0$
with respect to both $\xi_1'$ and $\xi_2'$, and 
$\xi_1'$ and $\xi_2'$ are both tight, from the Thurston--Bennequin
inequality we obtain 
\[ e(S', \xi_1')=  e(S', \xi_2')=0. \] 

In the complements of $N$ and $N_{\phi}$ we have $\xi_i|_{-(Y \setminus N)} = 
\xi_i'|_{-(Y_{\phi} \setminus N_{\phi})}$, therefore $e(S, \xi_i')= e(S, \xi_i)$. 
Since $\xi_1$ and $\xi_2$ are $C^0$--close to
$\F_1$ and $\F_2$ and $\partial S= \alpha_+ \cup \alpha_-$ is tangent to both $\xi_i$ and 
$\F_i$, we have $e(S, \xi_i)=e(S, \F_i)$ for $i=1,2$. 
From the additivity property of the relative Euler
class we obtain
\[ \langle c_1(\xi_i'), [\overline{S}] \rangle = e(S, \xi_i')+e(S', \xi_i')=e(S, \F_i), \]
therefore $\langle c_1(\xi_1'), [\overline{S}] \rangle \neq \langle c_1(\xi_2'),[\overline{S}] \rangle$.
\end{proof}

Lemma \ref{diverse} implies that the the  
$\spc$--structures ${\mathfrak s}_{\xi_1'}$ and ${\mathfrak s}_{\xi_2'}$
induced by $\xi_1'$ and $\xi_2'$ are not isomorphic, therefore 
$c^+(\xi_1')$ and $c^+(\xi_2')$ are linearly independent 
because $c^+(\xi_1') \in HF^+(Y, {\mathfrak s}_{\xi_1'})$ and $c^+(\xi_2') \in 
HF^+(Y, {\mathfrak s}_{\xi_2'})$. This implies 
that $c^+(\xi_1)$ and $c^+(\xi_2)$ are linearly independent too, 
therefore it proves Theorem \ref{nonnapapera}.

\section{Applications of Theorem \ref{nonnapapera}}
\subsection{sutured manifolds}
In order to apply Theorem \ref{nonnapapera} we need a way to construct
taut foliations in $3$--manifolds. This is provided by Gabai's sutured
manifold theory.
\begin{dfn}(\cite[Definition 2.6]{gabai:1})
A {\em sutured manifold} $(M, \gamma)$ is a compact oriented $3$--manifold $M$ 
together with a set $\gamma \subset \partial M$ of pairwise disjoint annuli $A(\gamma)$ and
tori $T(\gamma)$. Each component of $A(\gamma)$ is a tubular neighbourhood of an
oriented simple closed curve called {\em suture}. Finally every component
of $\partial M \setminus \gamma$ is oriented, and the orientations must be coherent with the 
orientations of the sutures.
\end{dfn} 
We define $R_+(\gamma)$ the subset of $R(\gamma)= \partial M \setminus \gamma$ where the orientation 
agrees with the orientation induced by $M$ on $\partial M$, and $R_-(\gamma)$ the 
subset of $\partial M \setminus \gamma$ where the two orientations disagree. We define also
$R(\gamma)=R_+(\gamma) \cup R_-(\gamma)$.

\begin{dfn}(\cite[Definition 2.10]{gabai:1})
A sutured manifold $(M, \gamma)$ is tight if $R(\gamma)$ minimises
the Thurston norm in $H_2(M, \gamma)$.
\end{dfn}

We will give the following definition only in the simpler case when 
no component of $\gamma$ is a torus, because this is the case we are 
interested in.
\begin{dfn}(\cite[Definition 3.1]{gabai:1} and 
\cite[Correction 0.3]{gabai:2})
Let $(M, \gamma)$ be a sutured manifold with $T(\gamma)= \emptyset$, and $S$ a properly 
embedded oriented surface in $M$ such that 
\begin{enumerate}
\item no component of $S$ is a disc with boundary in $R(\gamma)$
\item no component of $\partial S$ bounds a disc in $R(\gamma)$
\item for every component $\lambda$ of $S \cap \gamma$ one of the following holds:
\begin{enumerate}
\item $\gamma$ is a non-separating properly embedded arc in $\gamma$, or
\item $\lambda$ is a simple closed curve isotopic to a suture in $A(\gamma)$.
\end{enumerate}
\end{enumerate}
The $S$ defines a {\em sutured manifold decomposition}
\[ (M, \gamma) \stackrel{S} \leadsto (M', \gamma') \]
where $M'= M \setminus S$ and 
\begin{align*}
& \gamma' = (\gamma \cap M') \cup \nu (S_+ \cap R_-(\gamma)) \cup \nu (S_- \cap R_+(\gamma)), \\
& R_+(\gamma') = ((R_+(\gamma) \cap M') \cup S_+) \setminus {\rm int}(\gamma'), \\
& R_-(\gamma') = ((R_-(\gamma) \cap M') \cup S_-) \setminus {\rm int}(\gamma'),
\end{align*}
where $S_+$ and $S_-$ are the portions of $\partial M'$ corresponding to $S$
where the normal vector to $S$ points respectively out of or into 
$\partial M'$. 

A taut sutured manifold decomposition is a sutured manifold 
decomposition $(M, \gamma) \stackrel{S} \leadsto (M', \gamma')$ such that both
$(M, \gamma)$ and $(M', \gamma')$ are taut sutured manifolds.
\end{dfn}

\begin{dfn}(\cite[Definition 4.1]{gabai:1})
A {\em sutured manifold hierarchy} is a sequence of decompositions 
\[ (M_0, \gamma_0) \stackrel{S_1} \leadsto (M_1, \gamma_1) \leadsto \ldots \stackrel{S_n} \leadsto (M_n, \gamma_n) \]
where $(M_n, \gamma_n)=(R \times [0,1], \partial R \times [0,1])$ for some surface with boundary
$R$.
\end{dfn} 

The main results in sutured manifold theory are that for any
taut sutured manifold $(M, \gamma)$ there is a sutured manifold hierarchy
starting form $(M, \gamma)$ \cite[Theorem 4.2]{gabai:1}, and that we can 
construct a taut foliation on $(M, \gamma)$ such that $R(\gamma)$ is union of
leaves from a sutured manifold hierarchy starting form 
$(M, \gamma)$ \cite[Theorem 5.1]{gabai:1}. Thus sutured manifold theory 
translates the problem about the existence of taut foliations into a 
finite set of combinatorial data. The particular result we will use in 
our applications is the following.

\begin{thm}\label{costruzione}
Let $M$ be a closed, connected, irreducible, orientable $3$--manifold,
and let $\Sigma$ be a genus minimising connected surface representing
a non-trivial class in $H_2(M; \Q)$. Denote by $(M_1, \gamma_1)$ the taut 
sutured manifold where $M_1= M \setminus \Sigma$ and $\gamma_1= \emptyset$. If $g(\Sigma)>1$ and there 
is a properly embedded surface $S$ in $M_1$ yielding a taut sutured 
manifold decomposition, then $M$ admits a smooth taut foliation
${\mathcal F}$ such that:
\begin{enumerate}
\item $\Sigma$ is a closed leaf,
\item if $f$ is a representative of the germ of the holonomy map
      around a closed curve $\delta \subset \Sigma$, then 
\[ \frac{d^nf}{dt^n} (0)= \left \{ \begin{array}{ll} 
                                                1, \quad i=1, \\
                                		0, \quad i>1,
                                \end{array}
			\right. \]
\item $e({\mathcal F}, S)= \chi(S)$.
\end{enumerate} 
\end{thm} 
\begin{proof}
By \cite[Theorem 4.2]{gabai:1} $M$ admits a taut sutured manifold
hierarchy
\[ (M_1, \gamma_1) \stackrel{S} \leadsto (M_2, \gamma_2) \leadsto  \ldots \stackrel{S_n} \leadsto (M_n, \gamma_n), \]
then take the foliation ${\mathcal F}_1$ constructed from that 
hierarchy using the construction in \cite[Theorem 5.1]{gabai:1}. 
In particular ${\mathcal F}_1$ is smooth because $g(\Sigma)>1$ and $\partial M_1$ is
union of leaves. We obtain ${\mathcal F}$ by gluing the two components
of $\partial M_1$ together.

The smoothness of ${\mathcal F}$ along $\Sigma$ and part (2) come from the 
Induction Hypothesis (iii) in the proof of 
\cite[Theorem 5.1]{gabai:1}. Part (3) is a consequence of Case 2 of the
construction of the foliation in the
proof of \cite[Theorem 5.1]{gabai:1}; in fact ${\mathcal F}$ has a 
leaf which coincides with $S$ outside a small neighbourhood of $\Sigma$, and
spirals toward $\Sigma$ inside that neighbourhood.
\end{proof}

\begin{figure}\centering
\includegraphics[width=10cm]{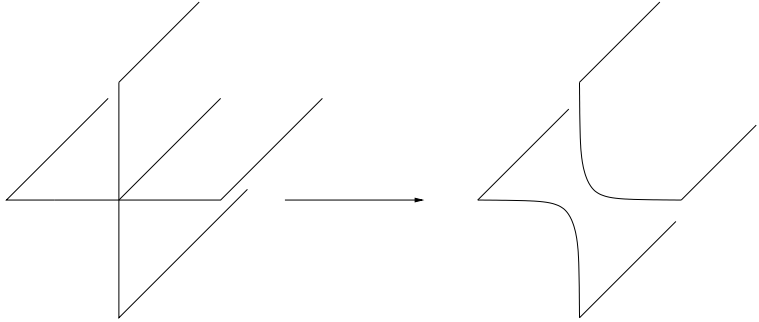}
\caption{Cut-and-paste surgery}
\label{cut-paste.fig}
\end{figure}

\begin{rem}\label{cut-paste}
A properly embedded surface $S$ in $M_1= M \setminus S$ gives a taut sutured 
manifold decomposition 
\[ (M_1, \gamma_1) \stackrel{S} \leadsto (M_2, \gamma_2) \] 
if for translates $\Sigma_+'$ and $\Sigma_-'$ of the boundary components of 
$\partial M \setminus \Sigma$ the surfaces $S+ \Sigma_+'$ and $S+ \Sigma_-'$ obtained by cut-and-paste 
surgery (see Figure \ref{cut-paste.fig}) are Thurston norm minimising 
in $H_2(M_1, \partial S)$. In fact $S+ \Sigma_+'$ and $S+ \Sigma_-'$ are isotopic to $R_+(\gamma_2)$
and to $R_-(\gamma_2)$ respectively, and minimising the Thurston norm in $M_1$
clearly implies minimising the Thurston norm in the smaller manifold
$M_2= M_1 \setminus S$.
\end{rem}

\subsection{Application to genus-one knots} \label{ss:nodi}
Let $K$ be a genus-one knot in $S^3$, and let $Y_K$ be the the 
$3$--manifold obtained as $0$--surgery on $K$. Let $T$ be a minimal 
genus Seifert surface for $K$ and let $\widehat{T}$ be the torus in 
$Y_K$ obtained by capping $T$ off. Denote $M_{\widehat{T}}= Y_K \setminus \widehat{T}$ 
and $\partial M_{\widehat{T}}= \widehat{T}_+ \cup \widehat{T}_-$, where $\widehat{T}_+$ is 
given the orientation induced by the orientation of $M_{\widehat{T}}$ by the
outward normal convention, and $\widehat{T}_-$ is given the opposite 
one. 

Let $\mu$ be a properly embedded curve in $M_{\widehat{T}}$ joining 
$\widehat{T}_+$ and $\widehat{T}_-$ which closes to the core of the 
surgery torus in $Y_K$, then $M_{\widehat{T}} \setminus \nu (\mu)$ is homeomorphic to 
$S^3 \setminus \nu (T)$. We can divide
$\partial (M_{\widehat{T}} \setminus \nu (\mu))$ in two pieces: $\partial_h (M_{\widehat{T}} \setminus \nu (\mu))=
\partial M_{\widehat{T}} \setminus \nu (\mu)$ called the {\em horizontal boundary}, and
$\partial_v (M_{\widehat{T}} \setminus \nu (\mu))= \partial \nu(\mu) \setminus \partial M_{\widehat{T}}$ called the {\em vertical boundary}.

Let $\alpha_+$ and $\beta_+$ be two simple closed curves in $\widehat{T}_+$ which
generate $H_1(\widehat{T}_+)$ and intersect transversally in a unique 
point.  Because the maps
\[ (\iota_{\pm})_* \colon H_1(\widehat{T}_{\pm}, \Z) \to H_1(M_{\widehat{T}}, \Z) \] 
induced by 
the inclusions are isomorphisms there are simple closed curves
$\alpha_-$ and $\beta_-$ in $\widehat{T}_-$ such that both $\alpha_+ \cup - \alpha_-$ and 
$\beta_+ \cup - \beta_-$ bound a surface in $M_{\widehat{T}}$. We may assume also that 
$\alpha_-$ and $\beta_-$ intersect transversally in a unique point.

Denote by ${\mathcal S}^+_n(\alpha)$ the set of the surfaces which are bounded 
by $\alpha_+ \cup - \alpha_-$ and which intersect $\mu$ transversally in exactly $n$
positive points and in no negative points, and by ${\mathcal S}^-_n(\alpha)$
 the set of the surfaces with the same property bounded by $- \alpha_+ \cup \alpha_-$. 
Let ${\mathcal S}^+_n(\beta)$ and ${\mathcal S}^-_n(\beta)$ be the same for the 
curves $\beta_+$ and $\beta_-$. Let $\kappa_n^+(\alpha)$ be the minimal genus of the
surfaces in ${\mathcal S}^+_n(\alpha)$ and  define $\kappa^-_n(\alpha)$, $\kappa_n^+(\beta)$, and 
$\kappa_n^-(\beta)$ in analogous ways. 

\begin{lemma}\label{asterix}
The sequences $\{ \kappa_n^+(\alpha) \}$, $\{ \kappa^-_n(\alpha) \}$, $\{ \kappa_n^+(\beta) \}$, and
$\{ \kappa_n^-(\beta) \}$ are non increasing.
\end{lemma}
\begin{proof}
We prove the lemma only for $\{ \kappa_n^+(\alpha) \}$ because the other
cases are similar. Let $S_n^+$ be a surface in ${\mathcal S}^+_n(\alpha)$
with genus $g(S_n^+)= \kappa_n^+(\alpha)$, and call $S_{n+1}^+$ the surface in 
${\mathcal S}^+_{n+1}(\alpha)$ constructed by cut-and-paste surgery
between $S_n^+$ and $\widehat{T}_+$. By definition $g(S_{n+1}^+) 
\geq \kappa_{n+1}^+(\alpha)$, and $g(S_{n+1}^+)= g(S_n^+)= \kappa_n^+(\alpha)$ because $\widehat{T}_+$ 
is a torus.
\end{proof}

\begin{lemma}\label{obelix}
If $K$ is not fibred, then for all $n \geq 0$ either $\kappa_n^+(\alpha) \neq 0$ and 
$\kappa^-_n(\alpha) \neq 0$, or $\kappa_n^+(\beta) \neq 0$ and $\kappa^-_n(\beta) \neq 0$.
\end{lemma}
\begin{proof}
Assume that there are annuli $A_{\alpha} \in {\mathcal S}_n^+(\alpha) \cup 
{\mathcal S}_n^-(\alpha)$ and $A_{\beta} \in {\mathcal S}_n^+(\beta) \cup {\mathcal S}_n^-(\beta)$.  
If we make $A_{\alpha}$ and $A_{\beta}$ transverse their 
intersection consists of one segment from $\alpha_+ \cap \beta_+$ to 
$\alpha_- \cap \beta_-$ and a number of homotopically trivial closed
curves. By standard arguments 
in three-dimensional topology we can isotope $A_{\alpha}$ and $A_{\beta}$
in order to get rid of the circles because $M_{\widehat{T}}$ is 
irreducible, therefore we can assume that $A_{\alpha} \cap A_{\beta}$ 
consists only of the segment. The boundary of 
$M_{\widehat{T}} \setminus (A_{\alpha} \cup A_{\beta})$ is homeomorphic to $S^2$, therefore  
$M_{\widehat{T}} \setminus (A_{\alpha} \cup A_{\beta})$ is homeomorphic to 
$(\widehat{T}_+ \setminus (\alpha_+ \cup \beta_+)) \times [0,1] \cong D^3$ because $M_{\widehat{T}}$ 
is irreducible. This proves that $M_{\widehat{T}}$ is homeomorphic to 
$\widehat{T}_+ \times [0,1]$. However, if $K$ is not fibred, then $Y_K$ is
not fibred either by \cite[Corollary 8.19]{gabai:3}, therefore 
$M_{\widehat{T}}$ is not a product.
\end{proof}

Assume without loss of generality that $\kappa_n^+(\alpha) \neq 0$ and 
$\kappa^-_n(\alpha) \neq 0$ for any $n \geq 0$. Lemma \ref{asterix} implies 
that the sequences $\{ \kappa_n^+(\alpha) \}$ and $\{ \kappa^-_n(\alpha) \}$ are definitively 
constant. Fix from now on a positive integer $m$ such that
$\kappa_{m+i}^+(\alpha)= \kappa_m^+(\alpha)$ and $\kappa_{m+i}^-(\alpha)= \kappa_m^-(\alpha)$ for all $i \geq 0$,
and choose surfaces $S_m^+ \in {\mathcal S}^+_m(\alpha)$ and 
$S_m^- \in {\mathcal S}^-_m(\alpha)$ such that $S_m^+$ has genus $\kappa_m^+(\alpha)$ and
 $S_m^-$ has genus $\kappa_m^-(\alpha)$.

\begin{lemma} \label{chirurgia}
Let $K_0$ and $K$ be knots in $S^3$, and denote by $M$ the $3$--manifold
obtained by gluing $S^3 \setminus \nu(K_0)$ to $S^3 \setminus \nu(K)$ via an 
orientation-reversing diffeomorphism $f \colon \partial (S^3 \setminus \nu(K_0)) \to \partial (S^3 \setminus \nu(K))$
mapping the meridian of $K_0$ to the meridian of $K$ and the longitude 
of $K_0$ to the longitude of $K$. Then $M$ is diffeomorphic to the 
$3$--manifold $Y_{K_0 \# K}$.
\end{lemma}
\begin{proof}
Glue a solid torus $S_1= S^1 \times [-\epsilon, \epsilon] \times [0,1]$ to $S^3 \setminus \nu(K_0) \sqcup S^3 \setminus \nu(K)$ 
so that $S^1 \times [- \epsilon , \epsilon] \times \{ 0 \}$ is glued to a neighbourhood of a 
meridian of $K_0$ and $S^1 \times [- \epsilon , \epsilon] \times \{ 1 \}$ is glued to a neighbourhood
of a meridian of $K$. The resulting manifold $M'$ is diffeomorphic to 
$S^3 \setminus \nu (K_0 \# K)$, so that if we glue a solid torus $S_2= S^1 \times D^2$ to 
$M'= S^3 \setminus \nu (K_0 \# K)$ mapping the meridian of $S_2$ to the longitude of 
$K_0 \# K$ we obtain $Y_{K_0 \# K}$.

Now we look at the gluing in the inverse order. First we glue $S_1$ to
$S_2$ so that $S^1 \times \{ - \epsilon , \epsilon \} \times [0,1] \subset \partial S_1$ is glued to disjoint
 neighbourhoods of two parallel longitudes of $S_2$, therefore the
resulting manifold is diffeomorphic to $T^2 \times [0,1]$. Then we glue
$S_1 \cup S_2 = T^2 \times [0,1]$ to $S^3 \setminus \nu(K_0) \sqcup S^3 \setminus \nu(K)$ and we obtain $M$. 
\end{proof}

\begin{proof}[Proof of Theorem \ref{principale}.]
In order to estimate the rank of $\widehat{HFK}(K,1)$ we 
cannot apply Theorem \ref{nonnapapera} directly; we have to increase 
the genus of $K$ artificially first. Let $K_0$ be any fibred knot with 
genus one, say the figure-eight knot. By 
\cite[Theorem 3]{gabai:murasugi1}
$K$ is fibred if and only if $K \# K_0$ is fibred. Also,
 the K\"unneth formula for connected sums Proposition \ref{kunneth},
\[ \widehat{HFK}(K \# K_0, 2) \cong \widehat{HFK}(K, 1) \otimes 
\widehat{HFK}(K_0, 1) \cong  \widehat{HFK}(K, 1), \]
and by Proposition \ref{relazione}
\[ HF^+(Y_{K \# K_0}, 1) \cong \widehat{HFK}(K \# K_0, 2). \]
Let $T$ and $T_0$ be genus minimising Seifert surfaces for 
$K$ and $K_0$ respectively, and call $\widehat{T}$ and
$\widehat{T}_0$ the corresponding capped-off surfaces in $Y_K$ 
and $Y_{K_0}$.
 
By Lemma \ref{chirurgia} the $3$--manifold $Y_{K \# K_0}$ is 
diffeomorphic to the union of 
$S^3 \setminus \nu (K)$ and $S^3 \setminus \nu (K_0)$ along the boundary via an
 identification $\partial (S^3 \setminus \nu (K)) \to \partial (S^3 \setminus \nu (K))$ mapping 
meridian to meridian and longitude to longitude. In $Y_{K \# K_0}$
the Seifert surfaces $T$ and $T_0$ are glued together to give a 
closed surface $\Sigma$ with $g(\Sigma)=2$ which minimises the genus in its 
homology class. Call $M_{\Sigma}= Y_{K \# K_0} \setminus \Sigma$. 
If we denote by $\mu_0$ a segment in $M_{\widehat{T}_0}$ which closes to
the core of the surgery torus in $Y_{K_0}$, we can see 
$M_{\Sigma}$ as $(M_{\widehat{T}} \setminus \nu(\mu)) \cup (M_{\widehat{T}_0} \setminus \nu(\mu_0))$ glued
together along their vertical boundary components.
 
From $S_m^+$ and $S_m^-$ we can construct
surfaces $\widehat{S}_+$ and $\widehat{S}_-$ in $M_{\Sigma}$ by gluing a 
copy of $T_0$ to each one of the $m$ components of 
$S_m^{\pm} \cap \partial_v (M_{\widehat{T}} \setminus \nu (\mu))$.
From an abstract point of view $\widehat{S}_+$ and $\widehat{S}_+$ are
obtained by performing a connected sum with a copy of $\widehat{T}_0$ 
at each of the $m$ intersection points between $S_m^+$ or $S_m^-$ and $\mu$,
therefore $g(\widehat{S}_{\pm})= \kappa_m^{\pm}(\alpha) +m$. 

Consider the taut sutured manifold $(M_{\Sigma}, \gamma)$ where $\gamma= \emptyset$. We claim 
that 
\[ (M_{\Sigma}, \gamma) \stackrel{\widehat{S}_+} \leadsto (M_{\Sigma} \setminus \widehat{S}_+, \gamma_+) \]
\[ (M_{\Sigma}, \gamma) \stackrel{\widehat{S}_-} \leadsto (M_{\Sigma} \setminus \widehat{S}_-, \gamma_-) \]
are taut sutured manifold decompositions. By Remark \ref{cut-paste}
this is equivalent  to proving that the surfaces $\widehat{S}_+ + \Sigma_+$, 
$\widehat{S}_+ + \Sigma_-$, $\widehat{S}_- + \Sigma_+$, and $\widehat{S}_- + \Sigma_-$ 
obtained by cut-and-paste surgery between $\widehat{S}_{\pm}$ and $\Sigma_{\pm}$ 
minimise the genus in their relative homology classes in 
$H_2(M_{\Sigma}, \alpha_+ \cup \alpha_-)$. 

We recall that $\widehat{T}_+$ and
$\Sigma_+$ are oriented by the outward normal convention, while 
$\widehat{T}_-$ and $\Sigma_-$ are oriented by the inward normal convention.
For this reason $\mu \cap \widehat{T}_+$ and $\mu \cap \widehat{T}_-$ consist
both of one single positive point. We consider only 
$\widehat{S}_+ + \Sigma_+$, the remaining cases being similar due to the 
above consideration.

Let $\widetilde{S} \subset M_{\Sigma}$ be a genus minimising surface in the 
same relative homology class as $\widehat{S}_+ + \Sigma_+$. We can see 
$\widetilde{S}$ as the union of two (possibly disconnected) 
properly embedded surfaces with boundary $S \subset M_{\widehat{T}} \setminus \nu(\mu)$ 
and $S_0 \subset M_{\widehat{T}_0}\ \setminus \nu(\mu_0)$, then 
$\chi(\widetilde{S})= \chi(S) + \chi(S_0)$. 
We can easily modify $\widetilde{S}$ without increasing its genus 
so that it intersects $\partial_v (M_{\widehat{T}} \setminus \nu(\mu))$ and 
$\partial_v (M_{\widehat{T}_0} \setminus \nu(\mu_0))$ in homotopically non trivial curves. 
The number of connected components of $\partial S_0= \partial S \cap \partial_v(M_{\widehat{T}} \setminus \nu(\mu))$ 
counted with sign is $m+1$.

Since $\chi(S) + \chi(S_0) = \chi(\widetilde{S})= 2- 2g(\widetilde{S})$, and
$\widetilde{S}$ is genus minimising, 
 both $S$ and $S_0$ maximise the Euler 
characteristic in their relative homology classes. 
$M_{\widehat{T}_0} \setminus \nu(\mu_0)$ is a product $T_0 \times [0,1]$, therefore 
$\chi(S_0)$ is equal to the negative of the number of components of 
$\partial S_0$ counted with sign, i.~e. $\chi(S_0)=-(m+1)$. We can modify $S_0$ 
without changing $\chi(S_0)$ so that it consists of some boundary 
parallel annuli and $m+1$ parallel copies of $T_0$, then we push 
the boundary parallel annuli into $M_{\widehat{T}} \setminus \nu(\mu)$, so that we have 
a new surface  $S' \subset M_{\widehat{T}} \setminus \nu(\mu)$ whose intersection with 
$\partial_v (M_{\widetilde{T}} \setminus \nu(\mu)$ consists of exactly $m+1$ positively oriented non 
trivial closed curves. If we glue discs to these curves we obtain
a surface $S_{m+1}^+ \in {\mathcal S}^+_{m+1}$ such that 
\[ g(\widetilde{S})= g(S_{m+1}^+ \# (m+1) \widehat{T}_0)= 
\kappa_{m+1}^+ + m+1. \]  
Since 
\[ g(\widehat{S}_+ + \Sigma_+)= g(\widehat{S}_+)+1= g(S_m^+ \# m \widehat{T}_0) 
+1= g(S_m^+)+m+1= \kappa_m^+ +m+1 \]
and $\kappa_{m+1}=\kappa_m$, we conclude that $g(\widehat{S}_+ + \Sigma_+)=
g(\widetilde{S})$, then $g(\widehat{S}_+ + \Sigma_+)$ minimises the 
genus in its relative homology class.

By Theorem \ref{costruzione} the taut sutured manifold decompositions
\begin{align*}
& (M_{\Sigma}, \gamma= \emptyset) \stackrel{\widehat{S}_+} \leadsto (M_{\Sigma} \setminus \widehat{S}_+, \gamma_+) \\
& (M_{\Sigma}, \gamma= \emptyset) \stackrel{\widehat{S}_-} \leadsto (M_{\Sigma} \setminus \widehat{S}_-, \gamma_-)
\end{align*}

provide taut smooth foliations $\F_+$ and $\F_-$ such that $\Sigma$ is a 
closed leaf for both so that, in particular,  
\[ \langle c_1(\F_+), [\Sigma] \rangle = \langle c_1(\F_-), [\Sigma] \rangle = \chi(\Sigma), \]
and moreover  
$e(\F_+, \widehat{S}_+)= \chi(\widehat{S}_+)$ and 
$e(\F_-, \widehat{S}_-)= \chi(\widehat{S}_-)$. 

Take $R \in {\mathcal S}_0^+(\alpha)$, then $-R \in {\mathcal S}_0^-(\alpha)$, therefore
 $[\widehat{S}_+]= [R]+ m[\Sigma]$ and $[\widehat{S}_-]= -[R]+ m[\Sigma]$ as 
relative homology classes in $H_2(M_{\Sigma}, \alpha_+ \cup \alpha_-)$, therefore 
\[ e(\F_+, \widehat{S}_+)= e(\F_+, R)+ m \chi(\Sigma) = e(\F_+, R)-2m \] 
and 
\[ e(\F_-, \widehat{S}_-)= e(\F_-, - R) + m \chi(\Sigma) =- e(\F_-, R)-2m. \]
This implies 
\begin{equation} \label{eq1}
e(\F_+, R)= \chi(\widehat{S}_+)+2m= \chi (S_m^+) = -2 \kappa_m^+(\alpha) 
\end{equation}
and 
\begin{equation} \label{eq2}
e(\F_-, R)= - \chi(\widehat{S}_-)- 2m= - \chi (S_m^-)= 2 \kappa_m^-(\alpha). 
\end{equation}
Recall that $\chi(S_m^{\pm})=- 2 \kappa_m^{\pm}(\alpha)$ because $S_m^+$ and $S_m^-$ have $2$ 
boundary components each. Equations \ref{eq1} and \ref{eq2} imply that 
$e(\F_+, R) \neq e(\F_-, R)$ because  $\kappa_m^{\pm}(\alpha)>0$, 
so we can apply Theorem \ref{nonnapapera}. 
\end{proof}

\begin{proof}[Proof of Corollary  \ref{nicecorollary1}.]
It is well known that the trefoil knots and the figure-eight
knot are the only three fibred knots with genus one, 
therefore by Theorem \ref{principale} if $\widehat{HFK}(K, 1)= \Z$
then $K$ is either one of the trefoil knots or the figure-eight knot.
These knots are alternating, therefore their knot Floer homology 
groups can be computed by using \cite[Theorem 1.3]{O-Sz:alternating},
 then the statement follows from the fact that these groups are 
distinct.
\end{proof}

\begin{proof}[Proof of Corollary  \ref{nicecorollary2}.]
The proof is immediate from from Corollary \ref{nicecorollary1} and 
from \cite[Theorem~1.6]{O-Sz:seifert} asserting that,
if surgery on a knot $K$ gives the Poincar\'e homology sphere, then the
knot Floer homology of $K$ is isomorphic to the knot Floer homology
of the left-handed trefoil knot.
\end{proof}

\bibliographystyle{plain}
\bibliography{contatto}
\end{document}